\documentclass[draft]{amsart}

\usepackage{amssymb}

\newtheorem{thm}{Theorem}[section]

\newtheorem{lem}[thm]{Lemma}

\theoremstyle{remark}
\newtheorem{rem}[thm]{Remark}

\theoremstyle{definition}

\numberwithin{equation}{section}
\numberwithin{thm}{section}

\begin{document}

%%%%%%%%%%%%%%%%%%%%%%%%%
% Subject classification 
%%%%%%%%%%%%%%%%%%%%%%%%%

% numbers and, optionally, one or more secondary classification numbers. 
% Use the following format:  "Primary 42B25. Secondary 42B60, 20E26"


%%%%%%%%%
% Title
%%%%%%%%%

% Title, in lower case, with no explicit linebreaks (\\).  If the title
% is too long to be used as a running head, add a short version of the
% title in brackets, as in \title[shorttitle]{fulltitle}.

\title{Diffusing Polygons and SLE($\kappa,\rho$)}

%%%%%%%%%%%%%%%%%%%%%%%%%%%%%%
% Author names and addresses 
%%%%%%%%%%%%%%%%%%%%%%%%%%%%%%

% Provide one separate \author{...} \address{...} \email{....} entry for each
% author, i.e., do not combine multiple authors.  Separate address lines by double
% slashes.  Do not attach footnotes to author  names. (For acknowledgements use
% the "\thanks" construct below.)
%

\author{Robert~O. Bauer}
\address{Department of Mathematics\\ University of Illinois at Urbana-Champaign\\ 1409 West Green Street \\ Urbana, IL 61801, USA}

\email{rbauer@math.uiuc.edu}

\author{Roland~M. Friedrich}
\address{Max-Planck Institut f\"ur Mathematik\\ D-53111
Bonn}

\email{rolandf@mpim-bonn.mpg.de}
%%%%%%%%%%%%%%%%%%%%
% Acknowledgements
%%%%%%%%%%%%%%%%%%%

% Use \thanks for acknowledgements as footnotes to the title page.  
% (Note that footnotes inside \author or \title macros are not
% allowed.)
%
% In case of multiple author papers, phrase the acknowledgement to 
% say "The first author was supported by ...  The second author was
% supported by ..."

\thanks{The research of the first author was supported by NSA grant H98230-04-1-0039.}

\thanks{The research of the second author was supported by
  NSF grant DMS-0111298.}

%%%%%%%%%%%%%
% Abstract 
%%%%%%%%%%%%%
%
% Abstracts should not contain macros (so that they can be processed independently
% of the paper.) Avoid displayed math and references in the abstract.

\begin{abstract}
We give a geometric derivation of $\text{SLE}(\kappa,\rho)$ in terms of
conformally invariant random growing compact subsets of polygons. The
parameters $\rho_j$ are related to the exterior angles of the polygons. We
also show that $\text{SLE}(\kappa,\rho)$ can be generated by a metric Brownian
motion, where metric and Brownian motion are coupled and the metric is a
pull-back metric of the Euclidean metric of an evolving polygon.      
\end{abstract}

\maketitle

\section{Introduction}

Stochastic Loewner evolution (or SLE) as introduced by Schramm in 
\cite{schramm:2000} describes random growing compacts in a simply connected 
planar domain $D$. Schramm discovered SLE  by considering discrete random simple
curves which satisfy (1) a Markovian-type property and whose scaling limit was
conjectured to be (2) conformally invariant. These two properties (plus a
reflection symmetry) render SLE canonical in the sense
that there exists only a one-parameter family of random non-self-crossing
curves $\gamma$ with these properties. They are all obtained by solving
L\"owner's equation \cite{loewner:1923} with a driving function given in terms
of Brownian motion. 
If $D$ is the upper half-plane
$\mathbb H$, and $\kappa\ge0$, consider for each $z\in\overline{\mathbb H}$ the ordinary
differential equation
\begin{equation}\label{E:CSLE}
\partial_t g_t(z)=\frac{2}{g_t(z)-W_t}, \quad g_0(z)=z,
\end{equation}
where $W_t=\sqrt{\kappa}B_t$, and $B_t$ is a 
one-dimensional standard Brownian motion. Let $T_z$ be the duration for which this
equation is well defined, i.e. $T_z=\sup\{t:\inf_{s\in[0,t]}|g_t(z)-W_t|>0\}$,
and set $K_t=\{z:T_z\le t\}$. Then it is easy to show that $g_t$ is a conformal
map from $\mathbb H\backslash K_t$ onto $\mathbb H$ with
$\lim_{z\to\infty}(g(z)-z)=0$. It can also be shown \cite{RS} that with
probability one the random
growing compact set $K_t$ is generated by a random non-self-crossing curve
$t\mapsto\gamma_t$ in the sense that $\mathbb H\backslash K_t$ is the unbounded
component of $\mathbb H\backslash\gamma[0,t]$. $\gamma$ is a random curve
connecting the boundary points $0$ and $\infty$ and is called {\em chordal}
$\text{SLE}_{\kappa}$ in $\mathbb H$ from $0$ to $\infty$. For an arbitrary domain $D$
and prime ends $z$ and $w$ chordal $\text{SLE}_{\kappa}$ in $D$ from $z$ to $w$ is 
defined via conformal invariance up to a time-change. I. e. if $f(\mathbb H)=D$, $f(0)=z$,
$f(\infty)=w$, then $f(\gamma_t)$ is a chordal
$\text{SLE}_{\kappa}$ in $D$ from $z$ to $w$. If $g$ is another conformal map
from $\mathbb H$ onto $D$ with $g(0)=z$, $g(\infty)=w$, then 
$\{g(\gamma_t),t\ge0\}$ has the same law as a time-change of $\{f(\gamma_t),t\ge0\}$. 

For calculations involving SLE conformal invariance is a powerful tool as it
is always permissible to choose the geometrically most convenient
configuration to do a given calculation. The solution depends only on the
conformal equivalence class, or the moduli, of the configuration. The
determination of certain hitting probabilities is reduced to solving
appropriate hypergeometric equations. In fact,
what one does is to track the evolution of these hitting probabilities as the
curve $\gamma$ grows, which comes down to tracking the evolution of the image
under the uniformizing map $g_t$ of the set $ \gamma$ is supposed to hit. For
example, if $\gamma$ is chordal $\text{SLE}_{\kappa}$ in the upper half-plane
from $0$ to $\infty$, $\kappa>4$ and $x,y>0$, then the probability that
$\gamma$ hits $(-\infty,-y)$ before $(x,\infty)$ depends only on the cross
ratio $-y/x$ and is given by 
\[
p=\frac{\Gamma(2-4a)}{\Gamma(2-2a)\Gamma(1-2a)}\left(\frac{y/x}{y/x+1}\right)^{1-2a}F(2a,
1-2a, 2-2a;\frac{y/x}{y/x+1}),
\]
for $\quad a=2/\kappa$, see \cite{lawlerbook}.
The calculation of $p$ uses the movement $x$ and $-y$ undergo under
the uniformizing map $g_t$, i.e. $t\mapsto x_t\equiv g_t(x)$, $-y_t\equiv
g_t(-y)$. We note that although $x_t$ and $y_t$ are coupled to $W_t$ via
\[
\partial_t x_t=\frac{2}{x_t-W_t},\quad\partial_t (-y_t)=\frac{2}{-y_t-W_t},
\]
there is no coupling of $W_t$ to $x_t$ or $-y_t$. If we do couple $W_t$ to
$x_t, -y_t$ via  
\[
dW_t=\sqrt{\kappa}dB_t+b(x_t,-y_t)dt,
\]
then the requirement that the random curve $\gamma$ that results from solving
L\"owner's equation for this $W_t$ be both, conformally invariant and satisfy a
Markovian-type property, forces the function $b$ to be homogenous of degree
$-1$, see \cite{BF:2005a}. The simplest such function is $b(x,y)=\rho_1/x+\rho_2/y$. Coupling with
this particular choice of drift $b$ leads to an example of
$\text{SLE}(\kappa,\rho)$, which we now define. 

Let $z_1< z_2<\dots< z_n$ be real numbers all distinct
from $0$. Consider the system of stochastic differential equations
\begin{align}\label{E:rho-diff}
        dW_t&=\sqrt{\kappa}\ dB_t
        +\sum_{k=1}^n\frac{\rho_k}{W_t-Z_t^k}\ dt,\notag\\
        dZ_t^k&=\frac{2}{Z_t^k-W_t}\ dt,\quad k=1\dots, n,
\end{align}
with $W_0=0, Z_0^1=z_1,\dots, Z_0^n=z_n$, and where $B_t$ is a one-dimensional
standard Brownian motion. The solution exists at least up to  some small
$t$. As above, let $g_t(z)$ be the solution to \eqref{E:CSLE}. Then the family
of conformal maps $g_t$ is called $\text{SLE}(\kappa,\rho)$ in the upper
half-plane from 
$(0,z_1,\dots, z_n)$ to $\infty$.
In this paper we will show that $\text{SLE}(\kappa,\rho)$ arises naturally
when one considers random growing compacts in polygons, i.e. that the
particular drift of $W_t$ in \eqref{E:rho-diff} can be derived from purely
geometric considerations. $\text{SLE}(\kappa,\rho)$, its properties and
relation to $\text{SLE}$ have been studied in several papers,
\cite{LSW:2003}, \cite{werner:2004}, \cite{dubedat:2003}.

%%%%%%%%%%%%%%%%%%%%%%%%%%%%%%%%%%%%%%%%%%%%%%%%%%%%%%%%%%%%%%

\section{Conformal Matters}\label{S:conformal}

\subsection{Schwarz-Christoffel Formula.} To begin with, let $D$ be a bounded simply connected domain whose boundary is a closed polygonal line without self-intersections. Let the consecutive vertices be $\zeta_1,\dots\zeta_n$ in positive cyclic order. The angle at $\zeta_k$ is given by the value of $\arg(\zeta_{k-1}-\zeta_k)/(\zeta_{k+1}-\zeta_k)$ between $0$ and $2\pi$ (we set $\zeta_{n+1}=\zeta_n$). Denote this angle by $\alpha_k\pi$, $0<\alpha_k<2$. We also introduce the outer angles $\beta_k\pi=(1-\alpha_k)\pi$, $-1<\beta_k<1$, and note that 
\begin{equation}\label{E:beta sum}
        \beta_1+\cdots+\beta_n=2.
\end{equation} 
The polygon is convex if and only if all $\beta_k>0$. We will call the pairs $(\zeta_k,\beta_k)$ the {\em corners} of the polygon.

Let $f$ be a conformal map from $D$ onto the upper half-plane $\mathbb H$ and let $z_k=f(\zeta_k)$. Assume that none of the $z_k$ equals $\infty$. For $z\in\mathbb H$ define the Schwarz-Christoffel mapping
\begin{equation}\label{E:sc}
        SC(z)=SC\left[\begin{array}{c | c}z_1, \dots, z_n &  z \\
        \beta_1, \dots,  \beta_n  & z^*\end{array}\right]=\int_{z^*}^z\prod_{k=1}^n(z-z_k)^{-\beta_k}\ dz,
\end{equation}
where the powers $(z-z_k)^{-\beta_k}$ denote analytic branches in $\mathbb H$. Note that 
\begin{equation}\label{E:SCprime}
        SC'(z)=\prod_{k=1}^n (z-z_k)^{-\beta_k},
\end{equation}
and
\begin{equation}\label{E:SC}
        \frac{SC''(z)}{SC'(z)}=-\sum_{k=1}^n\frac{\beta_k}{z-z_k}. 
\end{equation}
Then it is well known that for some constants $a,b\in \mathbb C$,
\[
        f^{-1}= a SC+b,
\]
see \cite{ahlfors:1966}. This result extends to the case when the polygon is allowed to have slits, i.e $\beta_k=-1$ for some $k$. Slits are counted as double edges of the boundary polygon, traversed in positive cyclic order. A vertex, when considered as a boundary point, may then occur multiple times, corresponding to different prime ends. Henceforth, a corner $(p_k,\beta_k)$ will always be a pair consisting of a prime end $p_k$ located at a vertex together with the exterior angle $\beta_k$ associated to the prime end $p_k$. The formula \eqref{E:sc} then remains unchanged if $z_k=f(p_k)$.

If $f(p_k)=\infty$ for one $k$, then \eqref{E:sc} needs to be adjusted by simply dropping the factor with exponent $-\beta_k$. 

Finally, formula \eqref{E:sc} continues to hold if $D$ is unbounded or one (or several) corners are at $\infty$, provided that the angles at $\infty$ are appropriately defined, see \cite{henrici:1974}. The interior angle at $\infty$ is chosen in $[-2\pi,0]$. For example, the angular region $\{z=r e^{i\varphi}:r>0, 0<\varphi<\alpha\}$ has a corner at infinity with angle $-\alpha$, while the infinite strip $\{0<\Im(z)<1\}$ has two corners at $\infty$ with angles $0$ (so $\beta=1$ there), and the slit plane $\mathbb C\backslash[0,\infty)$ has a corner at $\infty$ with angle $-2\pi$ ($\beta=3$). Then \eqref{E:beta sum} continues to hold. In the following, polygon will refer to these `generalized' polygons. 

%%%%%%%%%%%%%%%%%%%%%%%%%%%%%%%

\subsection{Marked domains, weights, and polygons}

We call a simply connected region $D$ together with a finite number of
distinct prime ends $p_1,\dots, p_n$ a {\em marked domain}. We call the set
$\{p_1,\dots, p_n\}$ the marking. A marked domain $(D,\{p_1,\dots,p_n\})$ is
{\em weighted} if to every prime end $p_k$ is associated a weight
$\beta_k\in\mathbb R$. We then call $\{(p_1,\beta_1),\dots,(p_n,\beta_n)\}$
the weighted marking. Two weighted marked domains 
\[
(D,\{(p_1,\beta_1),\dots,(p_n,\beta_n)\}),\quad\text{and}\quad(D',\{(q_1,\eta_1),\dots,
 (q_m,\eta_m)\})
\]
 are said to be conformally equivalent if  there is a conformal map $f$ from $D$ onto $D'$ such that  
\[
        \{(f(p_1),\beta_1),\dots, (f(p_n),\beta_n)\}
        =\{(q_1,\eta_1),\dots, (q_m,\eta_m)\}.
\]
In particular, this is only possible if $n=m$ and corresponding weights are equal.

If $D$ is a polygon it has a natural weighted marking given by its corners. We say two polygons $D$ and $D'$ are conformally equivalent if they are conformally equivalent as weighted marked domains using their natural weighted marking.

\begin{lem}
If $D$ is a polygon with corners $(p_1,\beta_1),\dots,(p_n,\beta_n)$ in positive cyclical order, and $\gamma$ is a Jordan arc contained in $D$ except for one endpoint which lies on the interior of a side $S$ of $D$, then there is a conformal map $f$ from $D\backslash\gamma$ onto a polygon $D'$ such that $(f(p_1),\beta_1),\dots, (f(p_n),\beta_n)$ are the corners of $D'$ in positive cyclical order. If we require $f(S\cup\gamma)=[0,1]$, then $f$ is unique.
\end{lem} 

\begin{proof}
Without loss of generality we may assume that $S\subset\mathbb R$. To begin, we also assume that $\gamma$ is differentiable, parameterized over the interval $[0,T]$ with $\gamma(0)\in S$. Then the existence of $f$ follows from an application of the Dirichlet principle: As $f$ is to preserve angles, $\arg f'$ is constant on $\partial D$. In fact if $f(S)\subset\mathbb R$, then $\arg f'=0$ on $\partial D$, and, for $t\in(0,T]$ we need $\arg f'=-\arg\dot{\gamma}(t)$. 

For general Jordan arcs $\gamma$, the result now follows by approximation.
\end{proof}

%%%%%%%%%%%%%%%%%%%%%%%%%%%%%%%

\subsection{Conformally invariant measures}

The following discussion is meant to motivate the ensuing construction. Let 
\[
D_{MW}\equiv(D,\{(p_1,\beta_1),\dots, (p_n,\beta_n)\})
\]
 be a weighted marked domain and $z$ and $w$ two prime ends of $D$ which are not part of the marking. Suppose that for each such triple $(D_{MW},z,w)$ we are given a probability measure $P_{D_{MW},z,w}$ on non-self-crossing curves in $D$ whose one endpoint is $z$ and whose other endpoint is $w$. Assume now that this family of probability measures is interrelated by {\em conformal invariance} and a {\em Markovian-type} property as follows:   

{\em Conformal invariance:} If $f:D\to f(D)$ is conformal, and $\gamma$ is a
random curve in $D$ from $z$ to $w$ with law  $P_{(D,\{(p_1,\beta_1),\dots,
  (p_n,\beta_n)\}),z,w}$, then $f(\gamma)$ is a random curve in $f(D)$ from
$f(z)$ to $f(w)$ with law 
\[
P_{(f(D),\{(f(p_1), \beta_1),\dots,(f(p_n),\beta_n)\}),f(z),f(w)}. 
\]
{\em Markovian-type property:} If $\gamma$ is a random curve in $D$ from $z$
to $w$ with law  $P_{(D,\{(p_1,\beta_1),\dots, (p_n,\beta_n)\}),z,w}$ and
$\gamma'$ is a sub-arc of $\gamma$ with endpoints $z$ and $z'(\in D)$, then
the conditional distribution of $\gamma$ given $\gamma'$ is 
\[P_{(D\backslash \gamma',\{(p_1,\beta_1),\dots, (p_n,\beta_n)\}),z',w}.
\]

Families of measures on curves in domains satisfying these two properties lead via Loewner's theory of slit mappings to random processes on the boundary of  appropriate domains. Conformal invariance and the Markovian-type property then imply that this motion together with the random motion of the conformal parameters---which classify the conformal equivalence class of the marked weighted domain---is a Markov process, see \cite{BF:2005a}. The question then centers on what Markov processes to consider. A particular choice leads to a well known  family of probability measures with the above properties---namely $\text{SLE}({\kappa},\rho)$---and our next goal is to derive the driving Markov process for $\text{SLE}({\kappa},\rho)$ from a purely geometric condition.

%%%%%%%%%%%%%%%%%%%%%%%%%%%%%%%%%%%%%%%%%%%%%%%%%%%%%%%%%%%%%%
 
\section{SLE$(\kappa,\rho)$ and Polygon motion}

Denote $\mathcal P$ the collection of all polygons $D$ such that 
\[
\{z:|z|<\epsilon,\Im(z)>0\}\subset D,\quad
\{z:|z|<\epsilon,\Im(z)=0\}\subset \partial D
\]
for some $\epsilon>0$. That is, $\mathcal P$ contains polygons with a side $l$
on the real axis, such that $0$ is in the interior of $l$ and such that the
interior of the polygon near $0$ lies in the upper half-plane.    
If $(p_1,\beta_1),\dots,(p_n,\beta_n)$ are the corners of a polygon 
$D\in\mathcal P$, we may assume that the corners are labeled such that
$0\in[p_1,p_2]\subset\mathbb R$. Then there is a map $SC:\mathbb H\to D$, 
\[
SC(z)=SC\left[\begin{array}{c | c}z_1, \dots, z_n &  z \\
        \beta_1, \dots,  \beta_n  & 0\end{array}\right]
\]
for some real points $z_1,\dots, z_n$ with $z_1<0<z_n$.
For $\kappa>0$ set
\begin{equation}\label{E:rho-mu}
        \rho_k=\frac{\kappa}{2}\beta_k,\quad k=1,\dots, n.
\end{equation}
In particular, $-\kappa/2\le\rho_k\le\kappa/2$.
Suppose that $(W_t, Z^1_t,\dots, Z_t^n)$ is a solution to \eqref{E:rho-diff}. 
 For $z$ in the upper half-plane, set
\[
        SC_t(z)=SC\left[\begin{array}{c | c}Z_t^1, \dots, Z_t^n & z \\
        \beta_1, \dots,  \beta_n & 0\end{array}\right].
\]
Then $z\mapsto SC_t(z)$ extends continuously to the real axis with the points $Z_t^k$ removed and is differentiable there as a function of $t$. In particular, if $W_s \neq Z_s^1,\dots, Z_s^n$ for $s\in[0,t]$, then we may define
\begin{equation}\label{E:def f}
        f_t(z)=SC_t(z)-\int_0^t(\partial_s SC_s)(W_s)\ ds.
\end{equation} 
Note that $f_t$ maps $\mathbb H$ into a polygon while the function $f$ in
Section \ref{S:conformal} mapped a polygon onto the upper half-plane.

Define the stopping time $\sigma$ by
\[
\sigma=\sup\{t:W_s, Z_s^1,\dots,Z_s^n\text{ are all distinct for }0\le s\le
t\}.
\]

\begin{lem}
The process $U_t\equiv f_t(W_t)$ is a martingale for $t<\sigma$. Furthermore, if 
\[
        A_t\equiv\kappa\int_0^t\left(SC'_s(W_s)\right)^2\ ds
\]
and $\tau(t)$ is defined by $A_{\tau(t)}=t$, then $t\mapsto U_{\tau(t)}$ is a standard Brownian motion.
\end{lem}

\begin{proof}
By an appropriate It\^o formula \cite{revuz.yor:1999},
\[
        dU_t=(\partial_t f_t)(W_t)\ dt+f_t'(W_t)\ dW_t
        +\frac{\kappa}{2} f_t''(W_t)\ dt
\]
Thus \eqref{E:def f}, \eqref{E:rho-diff}, and \eqref{E:SC} imply
\[
dU_t=\sqrt{\kappa} SC'_t(W_t)\ dB_t.
\]
By \eqref{E:SCprime}, $0<|SC'_t(W_t)|<\infty$ for $t<\sigma$, and the lemma follows.
\end{proof}

\begin{rem}
The motion of the corners of the polygon $f_t(\mathbb H)$ is
differentiable. Indeed, for $\epsilon$ small enough, 
\begin{align}
f_{t+\epsilon}(Z_{t+\epsilon}^l)-f_t(Z_t^l)%&=\int_{z^*}^{Z_{t+\epsilon}^l}
%\prod_{k=1}^n(z-Z_{t+\epsilon}^k)^{-\beta_k}\ dz-\int_{z^*}^{Z_{t}^l}
%\prod_{k=1}^n(z-Z_{t}^k)^{-\beta_k}\ dz\notag\\
&=\int_{\gamma_2}
\prod_{k=1}^n(z-Z_{t+\epsilon}^k)^{-\beta_k}\ dz+\int_{\gamma_3}
\prod_{k=1}^n(z-Z_{t+\epsilon}^k)^{-\beta_k}\ dz\notag\\
&\quad-\int_{\gamma_1}
\prod_{k=1}^n(z-Z_{t}^k)^{-\beta_k}\ dz-\int_t^{t+\epsilon}(\partial_s
SC_s)(W_s)\ ds,\notag
\end{align}
where $\gamma_1$ is a curve from $0$ to $Z_t^l$, $\gamma_2$ is a straight
line segment from $0$ to $Z_{t+\epsilon}^l-Z_t^l$, and $\gamma_3$ is
$\gamma_2$ shifted by $Z_{t+\epsilon}^l-Z_t^l$, connecting
$Z_{t+\epsilon}^l-Z_t^l$ to $Z_{t+\epsilon}^l$. If we parameterize
$\gamma_2$ by 
\[
s\in[0,1]\mapsto z(s)\equiv s(Z_{t+\epsilon}^l-Z_t^l),
\]
then 
\begin{align}
\lim_{\epsilon\to0}\frac{1}{\epsilon}&\int_{\gamma_2}
\prod_{k=1}^n(z-Z_{t+\epsilon}^k)^{-\beta_k}\ dz\notag\\
&=
\lim_{\epsilon\to0}\int_0^1\prod_{k=1}^n(s(Z_{t+\epsilon}^l-Z_t^l)-Z_{t+\epsilon}^k)^{-\beta_k}
\frac{Z_{t+\epsilon}^l-Z_t^l}{\epsilon}\ ds\notag\\
&=\frac{2}{Z_t^l-W_t}\prod_{k=1}^n(-Z_t^k)^{-\beta_k}.
\end{align}
Similarly, if $s\in[0,1]\mapsto z(s)$ is a smooth parameterization of
$\gamma_1$, and if 
\[
s\in[0,1]\mapsto \tilde{z}(s)\equiv z(s)+Z_{t+\epsilon}^l-Z_t^l
\]
parameterizes $\gamma_3$, then $\tilde{z}'(s)=z'(s)$ and
\begin{align}
\Delta&\equiv\int_{\gamma_3}
\prod_{k=1}^n(z-Z_{t+\epsilon}^k)^{-\beta_k}\ dz
-\int_{\gamma_1}
\prod_{k=1}^n(z-Z_{t}^k)^{-\beta_k}\ dz,\notag\\
&=\int_0^1\left(\prod_{k=1}^n(z(s)+Z_{t+\epsilon}^l-Z_t^l-Z_{t+\epsilon}^k)^{-\beta_k}-\prod_{k=1}^n(z(s)-Z_t^k)^{-\beta_k}\right)z'(s)\ ds.
\end{align}
Hence
\[
\lim_{\epsilon\to0}\frac{\Delta}{\epsilon}
=\frac{2(Z_t^l-Z_t^k)}{(Z_t^l-W_t)(Z_t^k-W_t)}\int_{0}^{Z_t^l}\prod_{k=1}^n
(z-Z_t^k)^{-\beta_k}\sum_{k\neq l}\frac{\beta_k}{z-Z_t^k}\ dz.
\]
\end{rem}

\begin{rem}
If we begin with an arbitrary $\text{SLE}(\kappa,\rho)$, i.e we begin with a
choice of $z_1,\dots, z_n$ and $\rho_1,\dots, \rho_n$, then the results of
this section continue to hold. In this case the Schwarz-Christoffel mapping $SC$ is no
longer guaranteed to be one-to-one. However, it still maps the intervals
$[z_k,z_{k+1}]$ onto straight line segments. By considering the Riemann
surface of the analytic function $SC$ we can still interpret the image
$SC(\mathbb H)$ as a polygon, albeit not a planar one. For example,
$\text{SLE}(2,(-1,-1))$, up to a normalization, leads to the map $z^3-3z$
which is easily understood in terms of a 3-fold cover, see \cite{ahlfors:1966}.
\end{rem}

For a polygon $D$, denote 
\[
(q,u)\in D\times\partial D\mapsto k_D(q,u)
\] the Poisson kernel of $D$. If $p\in\partial D$, denote $\partial_2
H_{D,p}(q,u)$ the analytic function in $q$ whose real part is $\partial_2
k_D(q,u)$ and which satisfies
\[
\lim_{q\to p}\partial_2 H_{D,p}(q,u)=0.  
\]

\begin{thm}[Loewner evolution in polygons]
Denote $K_t$ the hull of an $\text{SLE}_{\kappa}(\rho)$ in the upper half-plane
and $g_t:\mathbb H\backslash K_t\to\mathbb H$ the normalized uniformizing
map. Then $h_t\equiv f_t\circ g_t\circ f_0^{-1}:D\backslash
f_0(K_t)\to D_t$ satisfies
\begin{equation}\label{E:PL}
\partial_t\ln h_t'(z)=f_t(W_t)^2\partial_2 H_{D_t,f_t(\infty)}(h_t(z), f_t(W_t)).
\end{equation}
\end{thm}

\begin{proof}
Let $h_t=f_t\circ g_t\circ f_0^{-1}$. Then
\begin{align} 
h_t'(z)&=f_t'(g_t(f_0^{-1}(z)))g_t'(f_0^{-1}(z))(f_0^{-1})'(z) \notag\\
&=\frac{\prod_{k=1}^n(g_t(f_0^{-1}(z))-Z^k_t)^{-\beta_k} 
g_t'(f_0^{-1}(z))}{\prod_{k=1}^n(f_0^{-1}(z)-z_k)^{-\beta_k}}.\notag
\end{align}
Set $w=f_0^{-1}(z)$. As $\partial_t g_t'(z)=-2g_t'(z)/(g_t(z)-W_t)^2$, straightforward computation gives
\begin{align}
\partial_t h_t'(z)
&=\prod\left(\frac{g_t(w)-Z_t^k}{w-z_k}
\right)^{-\beta_k}g_t'(w)\notag\\
&\quad\times
\left[\sum_{l=1}^n\left(\frac{2}{g_t(w)-W_t}-\frac{2}{Z_t^l-W_t}
\right)\frac{-\beta_l}{g_t(w)-Z_t^l}
-\frac{2}{(g_t(w)-W_t)^2}\right]\notag\\
&=h_t'(z)\left[\frac{-2}{(g_t(w)-W_t)^2}
+\frac{2}{(g_t(w)-W_t}\sum_{l=1}^n\frac{\beta_l}{Z_t^l-W_t}\right].
\end{align}
Now, we note that
\[
H_{\text{Pol}_t}(q,u)=f_t'(f_t^{-1}(u))^{-1}H_{\mathbb H}(f_t^{-1}(q),f_t^{-1}(u)),
\]
whence, if $v_t=f_t^{-1}(u)$,
\begin{equation}
\partial_u H_{\text{Pol}_t}(q,u)
=f_t'(v_t)^{-2}\left[-\frac{f_t''(v_t)}{f_t'(v_t)}
H_{\mathbb H}(f_t^{-1}(q),v_t)+\partial_2 H_{\mathbb H}(f_t^{-1}(q), v_t)\right].
\end{equation}
Since $H_{\mathbb H}(z,w)=2/(z-w)$, the theorem
follows. 
\end{proof}

\section{SLE in variable background metric}

Instead of mapping $\text{SLE}(\kappa,\rho)$ into polygons we can also stay in
the upper half-plane and change the metric. Indeed, $f_t:\mathbb H\to D_t$ is
an immersion. If we endow $D_t$ with the Euclidean metric, then the metric
induced by $f_t$ on $\mathbb H$ is 
\[
g_{ij}=\delta_{ij}|f_t'(z)|^2,\quad i,j=1,2,
\]
where the indices 1 and 2 refer to the real and imaginary coordinate,
respectively. If $\Gamma=(\Gamma_{jk}^i)$ denotes the Levi-Civita connection
for this metric, then the (2-dimensional) Brownian motion $\tilde{W}$ for the
metric $(g_{ij})$ solves the stochastic differential equation
\[
d\tilde{W}_s^i=\sigma_j^i(\tilde{W}_s)\ dB_s^j-\frac{1}{2}
g^{kl}(\tilde{W}_s)\Gamma_{kl}^i(\tilde{W}_s)\ ds,
\]
see \cite{hsu:2002}. Here $g^{-1}=(g^{kl})$ is the inverse coefficient matrix
of $g$ and $\sigma$ is a square root of $g^{-1}$ (i.e. $\sigma
\sigma^T=g^{-1}$), and we observe the Einstein summation convention according
to which indices occurring once ``upstairs''  and once ``downstairs'' are to be
summed over. For our particular metric $g$ we find
\[
\Gamma_{11}^1=\Gamma_{22}^1=-\Re\left(\frac{f_t''}{f_t'}\right),
\]
see \cite{carmo:1992}. The boundary $\mathbb R=\partial\mathbb H$ is a 
one-dimensional sub-manifold of $\overline{\mathbb H}$. The metric $g$ on
$\mathbb H$ thus induces the metric $(f_t'(x))^2\ dx^2$ on $\mathbb R$. A
(one-dimensional) Brownian motion $W$ relative to this metric solves the
stochastic differential equation
\begin{equation}\label{E:metricBM}
dW_s=\frac{dB_s}{f_t'(W_s)}-\frac{1}{2(f_t'(W_s))^2}\sum_{j=1}^n\frac{\beta_j}{W_s-Z^k_t}\
ds.
\end{equation}
We now couple the metric to the Brownian motion $W$ via
\begin{equation}\label{E:couple}
dZ_t^k=\frac{2}{\kappa(f_t'(W_t))^2(Z^k_t-W_t)}\ dt,\quad k=1,\dots, n.
\end{equation}
Then, after a time-change, \eqref{E:metricBM} and \eqref{E:couple} become the
$\text{SLE}(\kappa,\rho)$-system \eqref{E:rho-diff} with the convention
$\rho_j=\kappa \beta_j/2$.

\end{document}